\documentclass[10pt,a4paper,reqno]{amsart}
\usepackage{amsfonts,amsthm,latexsym,amsmath,amssymb,amscd,amsmath,
epsf}
\usepackage{graphicx}

\newtheorem{Theorem}{Theorem}
\newtheorem{Lemma}{Lemma}

\newtheorem{Ex} {Example}

\newtheorem{theo+}           {Theorem}
\newtheorem{prop+}           {Proposition}
\newtheorem{coro+}           {Corollary}
\newtheorem{lemm+}           {Lemma}

\theoremstyle{definition}
\newtheorem{defi+}           {Definition}

\theoremstyle{remark}
\newtheorem{rema+}           {Remark}

\newenvironment{definition}{\begin{defi+}}{\end{defi+}}

\newcommand {\bCP} {\mathbb {CP}}

\newcommand{\bC}{\mathbb C}

\newcommand{\bZ}{\mathbb Z}

\newcommand{\RR}{\mathbb{R}}

\def \R{{\mathbb R}}

\def\newop#1{\expandafter\def\csname #1\endcsname{\mathop{\rm
#1}\nolimits}}

\newop{slc}
\newop{sign}
\newop{mdeg}

\begin{document}
          \numberwithin{equation}{section}

          \title[On limit sets for geodesics of meromorphic connections]
          {On limit sets for geodesics of meromorphic connections}

\author[D.~Novikov]{Dmitry Novikov}
\address{Faculty of Mathematics and Computer Science, Weizmann Institute of Science, Rehovot, 7610001 Israel}
\email{dmitry.novikov@weizmann.ac.il}

\author[B.~Shapiro]{Boris Shapiro}
\address{Department of Mathematics, Stockholm University, SE-106 91
Stockholm, Sweden}
\email{shapiro@math.su.se}

\author[G.~Tahar]{Guillaume Tahar}
\address{Faculty of Mathematics and Computer Science, Weizmann Institute of Science, Rehovot, 7610001 Israel}
\email{tahar.guillaume@weizmann.ac.il}

\date{\today}

\begin{abstract} Meromorphic connections on Riemann surfaces originate and are closely related to the classical theory of linear ordinary differential equations with meromorphic coefficients. Limiting behaviour of geodesics of such connections has been studied by e.g. Abate, Bianchi and Tovena  \cite{AB, AT} in relation with generalized Poincar\'{e}-Bendixson theorems. At present, it seems still  to be unknown whether some of the theoretically possible asymptotic behaviours of such geodesics really exist. In order to fill the gap, we use the branched affine structure induced by a Fuchsian meromorphic connection to present  several examples with geodesics having  infinitely many self-intersections and quite peculiar $\omega$-limit sets. 
\end{abstract}

\maketitle

\section{Introduction}

Meromorphic connections on complex analytic manifolds is both classical and modern area of mathematics interacting with a wide variety of topics including e.g., linear ordinary and partial differential equations in several  complex variables, theoretical and mathematical physics, differential geometry  and tensor calculus, representation and singularity theories. (For detailed expositions of the general theory of meromorphic connections an interested reader might consult \cite{HoTaTa, NoYa}). 

Motivated by the papers \cite{AB,AT, Ra}, in this short note  we discuss one aspect of meromorphic connections on Riemann surfaces, namely the asymptotic behavior of their self-intersecting geodesics. The structure of the paper is as follows. 
 
In \S~\ref{sec:conn}, we recall some general notions and facts related to meromorphic connections with a special emphasis on Fuchsian connections on tangent bundles of Riemann surfaces since they induce  branched affine structures which we shall use later.  
In \S~\ref{sec:hyperbolic}, we introduce the notion of a hyperbolic cylinder and use it to exhibit branched affine structures on $\bC P^1$ with geodesics accumulating on (possibly self-intersecting) limit cycles. 
In \S~\ref{sec:flat}, we focus on branched affine structures defined by  $k$-differentials. We provide several special examples of $k$-differentials on $\bC P^1$ and on the torus. These differentials have geodesics with both infinitely many self-intersections and minimal dynamics. In spite of this,  the presented geodesics  are still not dense in the ambient surfaces.

\medskip 
\noindent 
{\it Acknowledgments.}  We want to thank Dr.~Rakhimov for drawing our attention to this area. The first and the third author are supported by the Israel Science Foundation (grant No. 1167/17) and the European Research Council (ERC) under the European Union Horizon 2020 research and innovation program (grant agreement No. 802107). The second author wants to acknowledge the financial support of his research provided by the Swedish Research Council grant 2016-04416. \newline

\section{Meromorphic connections and branched affine structures}\label{sec:conn}

\subsection{Connections on vector bundles}
Let $\pi:E\to X$ be a holomorphic vector bundle on a complex manifold $X$.
A \emph{meromorphic connection} on $E$ is a $\bC$-linear operator $\nabla: \mathcal{M}_E\to \mathcal{M}_E\otimes\mathcal{M}^1_X$, where $\mathcal{M}_E$ is the sheaf of germs of meromorphic sections of $E$ and $\mathcal{M}^1_X$ is the sheaf of germs of meromorphic $1$-forms on $X$, satisfying the Leibniz rule
\begin{equation}\label{eq:leibnitz}
	\nabla(fs)=s\otimes\partial(f)+f\nabla(s),
\end{equation}
for all $s\in\mathcal{M}_E$ and $f\in \mathcal{M}_X$, the sheaf of meromorphic functions on $X$.

Let $s$ be a meromorphic section of $E$, $x\in X$ and $F_x=\pi^{-1}(x)$. The   $\nabla_x(s)$ is an $F_x$-valued one-form on $T_xX$ providing for every $\xi\in T_xX$ a vector $\nabla_x(s)(\xi)\in F_x$ considered at ``the derivative of $s$ at $s(x)$ in direction $\xi$".

\subsubsection*{Geometric point of view} The connection can be defined alternatively by a choice of a right inverse of $d_p\pi$ at every point $p\in E$, i.e. by the choice of lifting of vectors tangent to $X$ to points of $E$.  More exactly, let $p\in E$ and $x=\pi(p)\in X$. The inclusion $\iota: F_{x}=\pi^{-1}(x)\hookrightarrow E$ and the projection $\pi:E\to X$ canonically  induce the exact sequence 
\begin{equation}
	0\to  T_p F_{x}\stackrel{d_p\iota}{\to} T_pE\stackrel{d_p\pi}{\to} T_{x}X\to 0,
\end{equation}
however there is no canonical way to identify $T_pE$ with $ T_{x}X\times T_pF_{x}$.

A choice of a right inverse  $\sigma=\sigma(p): T_{x}X\to T_pE$ to $d_p\pi$ extends to the exact sequence 
\begin{equation*}
	0\to  T_{x}X\stackrel{\sigma}{\to} T_pE\stackrel{\tilde{\tau}}{\to}T_p F_{x} \to 0,
\end{equation*}
and defines an isomorphism 
$\left(d_p\pi, \tilde{\tau} \right):T_pE\cong T_{x}X\times T_pF_{x}$.
Combining $\tilde{\tau}$ with the canonical isomorphism $T_pF_{x}\cong F_{x}$ (here we use that $F_{x}$ is a linear space) we get the map  $\tau:T_pE\to F_x$   and the isomorphism 
\begin{equation*}
	\left(d_p\pi,\tau \right):T_pE\cong T_{x}X\times F_{x}.
\end{equation*}
The vectors in  $\operatorname{Im} d_p\iota, \operatorname{Im} \sigma$ are called vertical and horizontal vectors, respectively. 

Now, let  $s$ be a germ of a meromorphic section of $E$ at point $p=\left(x, s(x)\right)\in E$ and $\gamma:\left(\R,0\right)\to\left(X, x\right)$ be a germ of a curve at $x\in X$, with $\dot\gamma(0)=\xi\in T_{x}X$. Let
$\Xi=\frac{d \left(s(\gamma)\right)}{dt}(0)\in T_pE$ be  the tangent vector at $p$ tangent to the lifting $s(\gamma(t))$ of the curve $\gamma(t)$ to the section $s$,  $d_p\pi(\Xi)=\xi$. We define 
$$
\nabla_{x}(s)(\xi)=\tau(\Xi).
$$ 
If (as we always assume) the map $\sigma\left((x,v)\right)(\xi):F_x\times T_{x}X \to TE|_{F_{x}}$ depends $\mathbb{C}$-linearly on $v,\xi$  and  meromorphically on $x$ then thus defined  $F_x$-valued map $\nabla$ satisfies the Leibnitz rule \eqref{eq:leibnitz}.

\subsubsection*{Expressing a connection in local coordinates} 
A local trivialization $U\times F_x$, $x\in U\subset X$, of the bundle $E$ also induces an isomorphism $T_pE\cong T_xX\times F_x $, $p=(x,v)\in E$. Let $\{x_j\}$ be local coordinates on $U$, $\{e_k\}$ be a basis of $F_x$ and $s_k$ be the local sections $U \times \{e_k\}$. Then
\begin{equation}
	\sigma\left((x,s_k)\right)(\partial_{x_j})= \partial_{x_j}-\sum_i\Gamma_{jk}^i(x)s_i\in T_pE,
\end{equation}
where the \emph{Christoffel symbols} $ \Gamma^i_{jk}(x)$ are meromorphic functions on $U$. Thus 
$$
\nabla_x s_k(\partial_{x_j})=\tau(\partial_{x_j})=\partial_{x_j}-\sigma\left((x,s_k)\right)(\partial_{x_j})=\sum_i \Gamma^i_{jk}(x)s_i,$$
and 	therefore
\begin{equation}\label{eq:christoffel E}	
	\nabla_x s_k=	\sum_j\left(\sum_i\Gamma^i_{jk}(x)s_i\right)dx_j.
\end{equation}
Thus, for a section $s=\sum c_k(x)s_k$ we have
\begin{equation}\label{eq:nabla for s}
	\nabla_x (s)=\sum_k c_k(x) \nabla_x s_k+\sum_k s_k\otimes dc_k(x).
\end{equation}

\begin{Ex}\label{ex:standard connection} {\rm 
		Setting $\Gamma_{jk}^i=0$ in \eqref{eq:christoffel E} one obtains the \emph{standard connection} $\nabla_{st}$ on a local chart $U_i\times F_x$ and, in particular,    on $T\bC^n=\bC^n\times \bC^n$.  The horizontal sections of $\nabla_{st}$ are the constant sections $U_i\times \{v\}$ or the constant vector fields, respectively.
	}
\end{Ex}
We denote by $\Sigma=\Sigma(\nabla)$ the common \emph{singularity locus} of all $\Gamma^i_{jk}(x)$'s. If $\Sigma=\emptyset$ then the connection is called \emph{holomorphic}.
Tautologically, any meromorphic connection is holomorphic on $X\setminus\Sigma$.

\subsubsection{Horizontal paths and parallel transport}
A path $\tilde{\gamma}(t)=(\gamma(t), v(t))\subset E$ is called \emph{horizontal} if $\dot{\tilde{\gamma}}(t)$ is horizontal (i.e., $\nabla_{{\gamma}(t)}(\tilde{\gamma}(t)) \left(\dot{\tilde{\gamma}}(t)\right)=0$) for every $t$. 

Let $\gamma:[0,1]\to X\setminus\Sigma$ be a smooth path. 
\begin{Lemma}\label{lem:lifting horizontally}
	For any $v\in F_{\gamma(0)}$, there exists a horizontal lifting $\tilde{\gamma}:[0,1]\to E$ of $\gamma$ with $\tilde{\gamma}(0)=v$.
\end{Lemma}
\begin{proof}
	The induced vector bundle $\gamma^*E$ over $[0,1]$ inherits the induced connection $\gamma^*\nabla$.  The   condition $\gamma^*\nabla(s)=0$ for a section $s$ to be  horizontal  becomes, in a chosen trivialization, a system of 
	linear ordinary differential equations. A solution of this equation with the initial condition $\tilde{\gamma}(0)=v$ is the required horizontal lifting of $\gamma$.
\end{proof}

Horizontal liftings of $\gamma$ provide a global trivialization $E|_{\gamma(t)}\simeq [0,1]\times F_{\gamma(0)}$ of the restriction of  $E$ to the curve $\gamma(t)$,  and $\nabla$ becomes the standard connection in this trivialization.
In particular, this trivialization defines an isomorphism $\Gamma(\gamma)_0^1:  F_{\gamma(0)}\to F_{\gamma(1)}$  called the ``parallel transport of $F_{\gamma(0)}$ to $F_{\gamma(1)}$ along $\gamma$".

\subsubsection{Locally flat connections}

A holomorphic connection $\nabla$ is called \emph{locally flat} if there exist  local trivializations $U_i\times F$ of $E$ such that $\nabla $ becomes the standard connection in this chart. Equivalently, for a locally flat connection, there exists an   $\mathcal{O}(U_i)$-basis $s_i$  of the set of holomorphic sections of $E$ over $U_i$  called the \emph{basis of horizontal sections} such that $\nabla(s_i)=0$.  

More geometrically, local flatness  means that the parallel transport of $F_p$ to $F_{p'}$ along any curve $\gamma$ joining $p$ and $p'$ depends only on the homotopy  class of $\gamma$ (with fixed endpoints).
In particular, for a simply connected $X$, any locally flat holomorphic connection $\nabla$ on $E$  defines a global trivialization   $E\simeq X\times F$, with $\nabla$ being a standard connection in this trivialization.

Horizontal sections of a locally flat connection  on $E$ define a \emph{local system} on $X$, i.e., a locally constant sheaf associating to every open simply connected set $U\subset X$ the set of all  horizontal sections of $\nabla$ over $U$.
Vice versa, if for an atlas $\{U_i\times V\}$ of the vector bundle $E$, the transition functions $\psi_{ij}:U_i\cap U_j\to GL(V)$ are constant then the restriction of the  constant sections from $U_i$ and $U_j$ to $U_i\cap U_j$
coincide. Equivalently, the restrictions of $\nabla_{st}$ from $U_i\times F$ and $U_j\times F$  to $U_i\cap U_j\times F$ coincide, thus they glue together to a locally flat connection on $E$.

\smallskip 
For $\dim X>1$, the equation $\nabla(s_i)=0$ in local chart $U\times F$  becomes a system of partial differential equations on $U$ which can have no solutions. The integrability condition for this system  (i.e., the criterion of existence of a local basis of flat sections) is provided by the famous Frobenius theorem originally formulated in \cite{Fr}. The Frobenius condition can be equivalently reformulated as the condition of vanishing of the curvature tensor of the connection $\nabla$. 

For $\dim X=1$, i.e., in the case of Riemann surfaces, the system of partial differential equations $\nabla(s_i)=0$ for a horizontal section reduces to a system  of \emph{ordinary} differential equations. The latter is always solvable by the standard existence and uniqueness results for solutions of ordinary differential equations. In other words, for $\dim X=1$, a local basis of horizontal sections always exists and a holomorphic connection on a vector bundle over a Riemann surface is always locally flat.

\subsubsection{Tangent bundle of a Riemann surface}
From now on we concentrate on the case of the tangent bundle $TX$ of a Riemann surface $X$. 
%Natural interesting examples of holomorphic/meromorphic connections on $TX$ can be obtained as follows.  
%  Consider a polynomial  vector field $\mathcal{X}$ in $\bC^n$ such that the  Riemann surface $X$ (or a birationally equivalent algebraic curve)  is one of its global solutions. Then the restriction $\mathcal{X}_X$ becomes a holomorphic (or meromorphic) section of $TX$ thus defining a holomorphic (or meromorphic) connection on  $TX$, the one for which $\mathcal{X}_X$ is a horizontal section.

\subsubsection{Geodesics} Let  $\gamma(t)\subset X$ be (a germ of a) $C^2$-smooth curve on $X$. A choice of a holomorphic connection $\nabla$  on $TX$ implies,  by Lemma~\ref{lem:lifting horizontally}, the existence of  a horizontal lifting $\tilde{\gamma}(t)\subset TX$ of $\gamma(t)$ starting at $p=\left(\gamma(0),\dot{\gamma}(0)\right)$. 

Also, independently of the choice of $\nabla$, there is the
canonical lifting $j_1\gamma(t)=\left(\gamma(t), \dot{\gamma}(t)\right)\subset TX$ of $\gamma(t)$.  
\begin{definition}
	The curve $\gamma(t)\subset X$ is called \emph{geodesic} if these two liftings coincide.
\end{definition}
This means that the vector $\xi(t)=\dot{j_1\gamma}(t)\in T_{j_1\gamma(t)}TX$,  is horizontal for all $t$, i.e., $\nabla_{\gamma(t)}(\dot\gamma(t))(\dot\gamma(t))=0$.

Choosing a local trivialization and using \eqref{eq:nabla for s}, we see that this condition becomes an ordinary differential equation of second order. Thus, by the existence and uniqueness theorems for solutions of ordinary differential equations, for any point $x\in X$ and any $\xi\in T_xX$, there exists a geodesics $\gamma(t)$ with $\gamma(0)=x$ and $\dot\gamma(0)=\xi$.

Moreover, since the connection is holomorphic, the aforementioned ordinary differential equation of second order is real-analytic which implies that the geodesics are real-analytic as well.

\begin{Ex}\label{ex:geodesics of standard connection}
	For a standard connection on $T\mathbb{C}\cong\mathbb{C}\times\mathbb{C}$, the geodesics are exactly the straight lines $\gamma(t)=at+b$, $a,b\in\bC$. Indeed,  by Example~\ref{ex:standard connection} the horizontal sections are just the constant sections, i.e., $\dot{\gamma}(t)\equiv a$ is constant.  
\end{Ex}

\begin{Lemma}\label{lem:geodesics define the connection}
	Given  a Riemann surface $X$, the holomorphic connection on $TX$ is uniquely defined by  its geodesics.
\end{Lemma}
Indeed, for $\xi\in T_xX$ define $\sigma(\xi)\in T_{\xi}TX$ as the vector tangent to the lifting $j_1\gamma$ of the geodesic satisfying the initial conditions $\gamma(0)=x,\dot\gamma(0)=\xi$. As $\dim T_xX=1$, this completely defines the map $\sigma$ and therefore the connection.
\begin{rema+}
	For $\dim X>1$ the connection is uniquely defined by its geodesics if we assume that the connection is \emph{torsion free}, i.e., $\Gamma^i_{jk}=\Gamma^i_{kj}$ for all $i,j,k$.
\end{rema+}

\subsubsection*{Affine charts} Any open and simply connected chart $\phi:U\to \bC$, $U\subset X$, on a Riemann surface $X$   canonically extends to a chart $\Phi=(\phi, d\phi):TU\to T\bC\cong\bC\times\bC$ on the tangent bundle $TX$. By pulling back   the standard connection $\nabla_{st}$ on $T\bC$ one then defines a locally flat connection $\nabla=\Phi^*\nabla_{st}$ on $TU$. The  pull-backs $\left(d\phi\right)^{-1}(a\partial_z)$  of the horizontal sections of $\nabla_{st}$ are the horizontal sections of $\Phi^*\nabla$ and the pull-backs $\phi^{-1}(at+b)$ of geodesics $at+b$ of $\nabla_{st}$ are the geodesics of $\nabla$. 

We say that the chart $\phi$ is the  \emph{affine chart} of the connection $\nabla$.

\begin{Lemma}\label{lem:affine chart and geodesics}
	For any holomorphic connection $\nabla$ on $TX$ and any $x\in X$, there is an affine chart $\phi:U\to\bC$ for $\nabla$ defined in a neighborhood $U$ of $x$. 
\end{Lemma}
\begin{proof}
	Let $\gamma(t)$ with $\gamma(0)=x$, $\gamma'(0)\not=0$, be a germ of a geodesic of  $\nabla$. 
	By the implicit function theorem, there exists $\phi:(U,x)\to (\bC,0)$   such that  $\phi({\gamma}(t))\equiv t$. Then $\phi$ is an affine chart for $\nabla$: the curve $t$ is a geodesic for $\Phi_*\nabla$, so $\Phi_*\nabla=\nabla_{st}$ by Lemma~\ref{lem:geodesics define the connection} and Example~\ref{ex:geodesics of standard connection}.
\end{proof}

\begin{Lemma}\label{lem:transitions between affine charts}
	Let $\phi_i:U_i\to\bC$, $i=1,2$, be two affine charts on $X$. The transition map $\phi_{12}:\phi_2\left(U_1\cap U_2\right)\to\phi_1\left(U_1\cap U_2\right)$  has form  $\phi_{12}(z)=C_1z+C_2$ on any connected component of $U_1\cap U_2$.
\end{Lemma}
\begin{proof}  Indeed, $\phi$ maps geodesics of $\nabla_{st}$ to  geodesics of $\nabla_{st}$, i.e., $\phi_{12}(at+b)=a't+b'$ for some $a,b\in \bC$. \end{proof}

\begin{Lemma}\label{def:geodesic}
	A germ $\sigma:(\RR,0)\to (X,p)$
	of a  real curve at $p=\sigma(0)\in X$ is a {geodesic} of a holomorphic connection $\nabla$ on $TX$ if and only if for some (and thus for any) affine chart $\phi$ of $\nabla$ we have $\phi(\sigma(t))=at+b$ where $a,b\in\bC$.
\end{Lemma}

In local coordinates an affine chart looks as follows. Let $\nabla$ be a holomorphic connection on $TX$ and let $\tilde{\phi}:U\to\bC$ be a local chart on $X$ as above. Then $\tilde{\Phi}:TU\cong U\times\bC$ is a local chart on $TX$ and 
$$
\left(\tilde{\Phi}_*\nabla\right)(\partial_z)=\partial_z\otimes\eta
$$
for some one-form $\eta$ on $\phi(U)$. Let $s$ be a horizontal section of $\nabla$. Then    $\tilde{\Phi}(s)=f(z)\partial_z$  is a horizontal section of the connection $\tilde{\Phi}_*\nabla$ on $T\phi(U)$ and by \eqref{eq:leibnitz}
\begin{equation}\label{eq:horizontal section formula}
	0=\left(\tilde{\Phi}_*\nabla\right)(f \partial_z)= \partial_z\otimes (df+f\eta)\quad\text{which implies}\quad f=C\exp\left(-\int\eta\right).
\end{equation}
Finally, let ${\psi}:\phi(U)\to\bC$ be the map such that $\frac{d{\psi}}{dz}=\frac 1 f$, i.e.,
\begin{equation}\label{eq:affine chart}
	{\psi}=C^{-1}\int\exp\left(\int\eta\right)dz+C_2,
\end{equation}
and define $\phi:=\psi\circ\tilde{\phi}$.
Then the section $d({\phi})(s)=\partial_{{z}}$ is the horizontal section of $\Phi_*\nabla$. Thus $\Phi_*\nabla=\nabla_{st}$ and $\phi:U\to\bC$ is an affine  chart of $\nabla$.

\medskip

%As $\phi$ is a holomorphic chart, $\sigma$ is necessarily real analytic, and the above condition is equivalent to the requirement that the complexification $\sigma_\bC:(\bC,0)\to (X,p)$ is (a germ of) a horizontal section of $\nabla$.

%\subsection{Local metric}
%The standard Hermitian metric on $T\bC$ can be pushed forward to $TX$ by affine charts $\{U_\alpha\}$, giving a family of local metrics $\{g_\alpha\}$ \emph{adapted} to $\nabla$, i.e. such that 
%$$
%d \left(g_\alpha(R,T)\right)=\left(\nabla R,T\right)+\left(R,\bar{\nabla}T\right)
%$$ for any two germs of holomorphic vector fields $R,T$ on $X$. 

\subsection{Regular connections}
Let $\nabla$ be a meromorphic connection on a vector bundle $E\to X$ over a Riemann surface $X$, and let $p\in\Sigma$ be a pole of $\nabla$. 
%\begin{Def}
The connection $\nabla$ is called \emph{regular at $p$}  if in some chart $U\times V$ of $E$ the horizontal sections  $s$ have \emph{moderate  growth} at $p$, i.e., for any germ $\gamma(t)$ of a real-analytic curve at $p$, there exists $N>0$ such that $\|s(\gamma(t))\|=o\left(t^{-N}\right)$ as $t\to0$.
%\end{Def}

\subsubsection{Local theory} Here we consider only the tangent bundle $TX$. In this case the regularity property is equivalent to the condition that $\nabla$ is \emph{Fuchsian} at $p$, i.e., it has at most a simple pole at $p$.

Let $p$ be a simple pole of $\nabla$ and let $\alpha=\operatorname{res}_p\eta$ be the residue of $\eta$ at $p$. In a local  chart $z:U_p\to(\bC,0)$ near $p$ we have
$$
\nabla(\partial_z)=\partial_z\otimes\eta, \quad \eta=\frac\alpha z g(z)\,dz,
$$ 
where $g(z)$ is a holomorphic function  such that $g(0)=1$. 
\begin{Lemma}
	After a suitable biholomorphic change of variable $z$ we can assume that $g\equiv1$ and    $\eta=\frac\alpha z dz$.
\end{Lemma} 
\begin{proof}
	By integrating  $\frac\alpha z g(z)\,dz=\frac\alpha w \,dw$ we get 
	$$\log w=\log z+\int g_1(z)dz,$$ where $g_1(z)=\frac{g(z)-1}{z}$ is holomorphic at $0$.
	Thus $w=z\exp\left(\int g_1(z)dz\right)=z+\dots$ is the required change of variables.
	\end{proof}
Therefore by \eqref{eq:horizontal section formula} the horizontal sections of $\nabla$ in this chart are given by  
$Cz^{-\alpha}\partial_z$  and, by \eqref{eq:affine chart}, the affine charts are 
\begin{align}\label{eq:affine near pole}
	w=\phi(z)=&C_1z^{-\alpha+1}+C_2 \quad \text{if} \quad \alpha\not=1, \quad or\\
	w=\phi(z)=&C_1\log z+C_2 \quad \text{if}\quad\alpha=1.
\end{align}  For $\alpha\not\in\bZ$ and for $\alpha=-1$, these maps are branching at $p$ which implies that there is no uni-valued affine chart in a neighborhood of $p$.  However, one can always cover the punctured neighborhood $D_\circ$ of $P$ by open  simply connected domains and define the affine charts for them. By Lemma~~\ref{lem:transitions between affine charts}, in the connected components of the intersections of these domains, the affine charts will be related by affine transformations.

For example, cover $D_\circ$ by the union of two sectors $S_-=\{0<|z|<\epsilon, |\arg z|<\pi\}$  and $S_+=\{0<|z|<\epsilon, 0<\arg z<2\pi\}$. The intersection $S_-\cap S_+=U_+\cup U_-$ consists of the union of two sectors, namely,  $U_-=\{0<|z|<\epsilon, \pi<\arg z<2\pi\}$ and $U_+=\{0<|z|<\epsilon, 0<\arg z<\pi\}$. 
Let $\phi_+, \phi_-$ be uni-valued branches of the affine chart \eqref{eq:affine near pole} in $S_+, S_-$ respectively. 
If $\phi_-= \phi_+$ in $U_+$ then on $U_-$ these two different holomorphic continuations of the same affine chart on $U_-$ are related by the linear transition map $\phi_{+-}(w)=\phi_+\circ\phi_-^{-1}(w)$, where
\begin{align}\label{eq:affine transition}
	\phi_{+-}(w)	=&e^{2\pi i\alpha}w+C_2\left(1-e^{2\pi i\alpha}\right),\quad\alpha\not=-1,\quad \textrm{or}\\
	\phi_{+-}(w)	=& w+2\pi i C_1, \quad\alpha=-1.
\end{align}

% $\phi(z)=z^{-\alpha+1}+C$, $\alpha\not\in \mathbb{Z}$, has univalued branches $\phi_-, \phi_+$ in the

The multiplicator  $|e^{2\pi i\alpha}|=e^{-2\pi Im(\alpha)}$  is called the \emph{dilation ratio} of $\nabla$ at $p$.

\subsubsection{Global theory}
Let $X$ be a  Riemann surface and $\nabla$ be a holomorphic connection on $TX$. 
The parallel transport defines the monodromy representation  $M:\pi_1 \left(X , p_0\right)\to GL(T_{p_0}X)\cong \bC^*$,  the monodromy of $\nabla$ along $\gamma$.

\smallskip 
Let $\{U_i\}$ be a covering of $X$ by affine charts $\phi_i$, and assume that $U_i\cap U_j$ are connected. The transition maps 
$$
\phi_{ij}:\phi_j(U_i\cap U_j)\to\phi_i(U_i\cap U_j), \quad \phi_{ij}(z)=\phi_i\circ\phi_j^{-1}(z)
$$
are affine by Lemma~\ref{lem:transitions between affine charts}. 

\begin{Lemma}
	One can choose  affine charts $\psi_i$ such that all transition maps $\psi_{ij}$ will be of the form
	$$
	\psi_{ij}(z)=a_{ij}z+b_{ij}, \quad a_{ij}\in M\left(\pi_1 \left(X , p_0\right)\right)\subset\bC^*.
	$$
\end{Lemma}

\begin{proof}
	Indeed, for a simply connected $X$, the bundle $TX$ can be globally trivialized by parallel transport, $TX\equiv X\times \bC$, with the section $s=X\times \{1\}$ being the horizontal section of $\nabla$. Let $\phi_i$ be an affine chart on $X$. The image $d\phi_i(s)$ is a constant vector field $a_i\partial_z$. Replacing $\phi_i$ by $a_i^{-1}\phi_i$, one can assume that $a_i=1$ for all $i$. Then the transition maps are given by  $\psi_{ij}(z)=z+b_{ij}$ since the  differentials $d\psi_i, d\psi_j$ agree on  every fiber $T_pX$, $p\in U_i\cap U_j$. 
	
	For a non-simply connected $X$, consider its universal cover $\pi:\tilde X\to X$. The lifts $U_{i,\iota}$ of $U_i$ define a cover of $\tilde{X}$ with the same transition maps $\phi_{ij}$.
	Repeating the previous step for the lifting of $\nabla$ to the $T\tilde{X}$, choose new affine charts $\psi_{i,\iota}$ of $U_{i,\iota}$ and define the affine charts $\psi_i:=\psi_{i,\iota}\circ\left(\pi_{U_i}\right)^{-1}$  for some arbitrary $\iota$.
	By definition, the maps $d\psi_{i,\iota}\circ d\psi_{i,\iota'}^{-1}$ present the monodromy of $\nabla$ along some path and the claim follows.
	
\end{proof}

\begin{defi+}
	An affine structure on a Riemann surface $X$ branched at a discrete subset $\Sigma\subset X$ is a finite holomorphic atlas on a $X\setminus\Sigma$ whose charts $U_i\subset X$ have piecewise real-analytic boundary and  transition maps are affine.
	We wil call a Riemann surface with a branched  affine surface a \emph{branched affine surface}.
	
	The affine structure is called \emph{Fuchsian} if for every $p\in\Sigma$ and any affine chart $\phi:U\to\bC$, $p\in \bar{U}$, there exists a (finite or infinite) limit $\lim_{U\ni z\to p}\phi(z)$.
\end{defi+}

\smallskip
Summing up, we obtain the following. 
\begin{Theorem}\label{thm:structure from connection}
	{\rm (i)} 	Affine charts of any Fuchsian connection $\nabla$ on $TX$ define a Fuchsian affine structure on $X$ branched at the polar locus $\Sigma$. 
	
	\noindent	
	{\rm (ii)} 	The multiplicators of the affine transition maps belong to the monodromy group $M\left(\pi_1 \left(X\setminus\Sigma, p_0\right)\right)\subset\bC^*$. 
	
	\noindent	
	{\rm (iii)} 	In particular, for $X=\bCP^1$, the multiplicators  belongs to the multiplicative subgroup of $\bC^*$ generated by $e^{2\pi i\alpha_j}$, where $\alpha_j$ are the residues of $\nabla$ at its poles.
\end{Theorem}

The opposite holds as well.

\begin{Theorem}\label{thm:connection from structure}
	An affine structure on a Riemann surface $X$ branched at a discrete set $\Sigma\subset X$ defines a connection $\nabla$ on $T(X\setminus\Sigma)$, with charts of the structure being the affine charts of $\nabla$.
	
	If the structure is Fuchsian, then the connection  $\nabla$ is Fuchsian as well.
\end{Theorem}
\begin{proof}
	Using charts $\phi_i:U_i\subset X\setminus\Sigma\to\bC$ we can pull-back the standard connection $\nabla_{st}$ on $T\left(\phi_i(U_i)\right)$ to define a connection $\nabla_i$ on $TU_i$. Since the transition maps $\phi_{ij}$ are affine, the geodesics of $\nabla_i$ and $\nabla_j$ on the intersections $U_i\cap U_j$ coincide: they both are sent to the straight lines $at+b$ by $\phi_i, \phi_j$, respectively. By Lemma~\ref{lem:geodesics define the connection}, the $\nabla_i$ and $\nabla_j$  then coincide on $T(U_i\cap U_j)$, so the connections $\nabla_i$ glue to  a connection $\nabla$ on  $T(X\setminus\Sigma)$.
	
	Assume now that the structure is regular and let $D_\circ$ be a punctured disc with center at $p\in\Sigma$ covered by some affine charts $\phi_i:U_i\to\bC$, i.e., $D_\circ=\cup_{i=1,\dots,n}U_i$.  Shrinking $D_\circ$ and $U_i$ we can assume that $U_i$ are sectors of $D_\circ$ ordered counterclockwise. 
	\begin{Lemma}
		The functions $\phi_i(z)$ have moderate growth at $p$ (here $z$ is a local parameter at $p$).
	\end{Lemma}
\begin{proof}
	Postcomposing  $\phi_i(z)$  with affine maps, we can ensure that the transition maps $\phi_{i,i+1}$ are identity for $i=1,\dots,n-1$ and $\phi_{n1}(w)=e^{2\pi i \alpha}w+C$. Then  the functions $\left(z^{\alpha}+C'\right)\phi_i(z)$ (or  $C'\log z+\phi_i(z)$ if $\alpha\in\bZ$) coincide on $U_i\cap U_j$ for a suitable $C'$ and are therefore restrictions of a function $\phi$ holomorphic in $D_\circ$. 
	
	Clearly, the limits $\lim_{U_i\ni z\to p}\phi(z)$ are all either simultaneously infinite or simultaneously finite. Thus either $\phi(z)$ or $1/\phi(z)$ is bounded as $z\to 0$ and therefore holomorphic at $p$. Then $\phi(z)$ (and therefore $\phi_i(z)$ as well) has moderate growth at $p$.
\end{proof}
The second claim  now follows easily: the pull-backs of constant sections of $T\bC$ (horizontal sections of $\nabla_{st}$) by maps with moderate growth have moderate growth.
\end{proof}

\subsection{Branching points}

Poles of a  meromorphic Fuchsian connection are the branching points of the corresponding affine structure. 

\smallskip 
Computations of \S~2.1.1 can be interpreted in the following way. 
Given a singularity $p$, the geometric model of the branched affine structure around $p$ is represented by an angular sector with the total angle $\theta$ where the extremal rays are identified by a homothety with ratio $\lambda$, see Lemma~\ref{lem:transitions between affine charts}. A natural local coordinate on this sector is of the form $f(z)=z^{1+\alpha}$ where $Re(\alpha)=\frac{\theta}{2\pi}-1$ and $Im(\alpha)=-\frac{\log \lambda}{2\pi}$. 
 The affine Schwarzian (invariant under postcompositions with affine transformations) is given by $\frac{f''}{f'}=\frac{\alpha}{z}$ where $-\alpha$ is the residue of the Fuchsian meromorphic connection at $p$.\newline 

Just like the fact that the conical angle can be interpreted as the discrete analog of curvature in the Gauss-Bonnet formula, the dilation ratio around a singularity can be thought of as a discrete imaginary curvature.

\subsection{Geodesics}

A curve on a branched affine surface is called a \textit{geodesic} if it is locally conjugated to a straight line in any affine chart, see Lemma~\ref{def:geodesic}. This class of curves coincides with the class of geodesics defined by the Fuchsian meromorphic connection on the tangent bundle corresponding to the branched affine structure under consideration.\newline

A \textit{regular geodesic} is a geodesic that does not meet any branching point of the affine structure.\newline

A \textit{saddle connection} (not to be confused with other connections we introduced earlier!) is a (possibly self-intersecting) geodesic segment whose endpoints are branching points of the branched affine structure or, equivalently, the singularities of the meromorphic connection.

\smallskip
Let us finally introduce the most crucial notion of this paper, namely, the \emph{$\omega$-limit set} of a geodesic. Informally, it is the set to which a geodesic accumulates. More exactly, the definition is as follows. 

\begin{defi+} Let $\gamma: \mathbb{R}^{+} \rightarrow X$ be a geodesic curve on a branched affine surface (or, equivalently, on a Riemann surface with a meromorphic connection). We define the $\omega$-limit set of $\gamma$ as $\bigcap\limits_{t \in \mathbb{R}^{+}} \overline{\gamma([t,+\infty[)}$.
\end{defi+}

\section{Hyperbolic cylinders and examples}
\label{sec:hyperbolic}

\subsection{Cylinders}

\begin{Lemma}\label{lemma:regularcylinder}
 Every regular closed geodesic without self-intersections on a branched affine surface  belongs to a $1$-parameter family of homotopy equivalent closed regular geodesics.
\end{Lemma}

\begin{proof}
We consider a neighborhood $V$ of a regular closed geodesic $\gamma$ on a branched affine surface. Without loss of generality, we can assume $V$ is an open topological cylinder without  branching points. The monodromy along $\gamma$ preserves its local direction. Therefore, the multiplicator of the affine map induced by the monodromy along $\gamma$ is a real positive number $a$. If $a=1$, there is a neighborood of $\gamma$ inside $V$ that is affinely equivalent to a parallelogram of $\mathbb{C}$ where a pair of sides is identified by a translation. If $a \neq 1$, then there is a neighborhood of $\gamma$ inside $V$ that is affinely equivalent to a trapezoid of $\mathbb{C}$ where a pair of sides is identified by a homothety. In both cases, $\gamma$ belongs to a continuous family of homotopic geodesics. 
\end{proof}

For any regular closed geodesic, we refer to the maximal family of geodesics it embeds into as a \textit{cylinder}.

\subsection{Hyperbolic cylinders and limit cycles}

A \textit{Hopf torus} is an annulus in the complex plane $\mathbb{C}$ whose inner and outer boundaries are identified by a homothety, see Subsection 2.2 in \cite{DuFoGh}.

A \textit{hyperbolic cylinder} of angle $\theta$ and dilation ratio $\lambda$ is an angular portion of a Hopf torus whose angle equals $\theta$ and whose dilation coefficient equals $\lambda$. Since this surface has an atlas with transition maps of the form $z \rightarrow az+b$ where $a \in \mathbb{R}_{+}^{\ast}$ this topological cylinder has an affine structure. Besides, since homotheties preserve directions, the foliation of a hyperbolic cylinder by the straight lines having the same slope is globally well-defined.

\begin{figure}

\begin{center}
\includegraphics[scale=0.7]{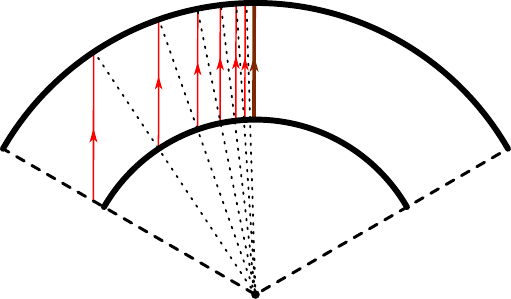} 
\end{center}

\vskip 0.5cm

\caption {Hyperbolic cylinder}
\label{fig:ill1}
\end{figure}

A hyperbolic cylinder of angle $\theta < \pi$ is affinely equivalent to a trapezoid whose identified sides are segments of $\mathbb{C}$ between the endpoints of the two circular arcs. This implies that every regular closed geodesic whose monodromy coincides with a nontrivial homothety belongs to a hyperbolic cylinder (see Lemma~\ref{lemma:regularcylinder}).\newline

In the case of a hyperbolic cylinder, the closed geodesics are the rays of the annulus. Any oriented geodesic entering such a cylinder with the same direction as a closed geodesic accumulates inside the cylinder on an attracting limit cycle, see Figure 1. Indeed, since directions are well-defined in a hyperbolic cylinder, a geodesic cannot cross a ray of the annulus that belongs to the same direction.\newline

In the latter case, for a hyperbolic cylinder covering the interval $[\alpha,\beta] \subset \mathbb{S}^{1}$ of directions and bounded by saddle connections, a geodesic in direction $\beta$ entering  the cylinder by crossing the boundary of direction $\alpha$ will accumulate onto the saddle connection forming the boundary of direction $\beta$. We will need the following result. 

\begin{Theorem}[Theorem 0.1 of \cite{AT}]
For an infinite geodesic $\gamma$ on $\mathbb{CP}^{1}$ without self-intersections, the $\omega$-limit set of $\gamma$ is either

\noindent
{\rm(i)} a closed geodesic

or

\noindent
{\rm (ii)} a singular limit cycle formed by saddle connections.
\end{Theorem}

Example 8.1 of \cite{AT} contains a numerical experiment illustrating Case (i) while there is no example in loc. cit. to illustrate  Case (ii). However, geodesics entering a hyperbolic cylinder confirm the existence of both cases, see Figure 1. (Theorem~2 proves that the branched affine surface containing a hyperbolic cylinder has a well-defined Fuchsian meromorphic connection where the geodesics are the same.)\newline

\subsection{Self-intersecting cylinder}

In our next example, we consider a branched affine structure on $\mathbb{CP}^{1}$ corresponding to a meromorphic Fuchsian connection such that the residues at its poles have their real parts belonging to $\frac{1}{4}\mathbb{Z}$. Following Subsection 2.3, this means that directions in the affine structure are well-defined up to rotations of order four.\newline

 Such an affine structure admits  a ramified covering of degree at most four on which all directions are well-defined. This canonical cover (discussed in Section 6 of \cite{Ra}) is a generalization of the canonical cover of Subsection 4.1. Then, every regular closed geodesic, even if it itself has self-intersections, will have no self-intersection on the finite cover since on this cover  it belongs to a cylinder, see Lemma~\ref{lemma:regularcylinder}. The projection of this cylinder to the original surface is a \textit{self-intersecting cylinder}.\newline

\smallskip

 \begin{figure}

\begin{center}
\includegraphics[scale=0.55]{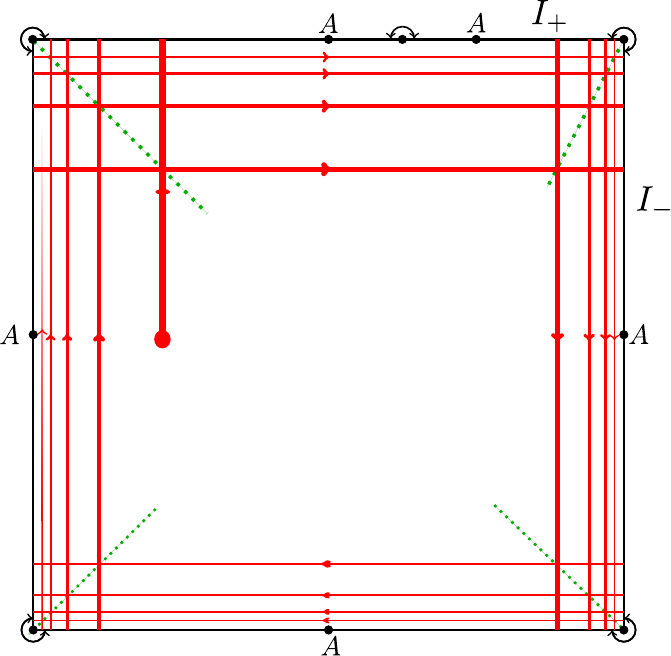} 
\end{center}

\vskip 0.5cm

\caption {Self-intersecting cylinder: the segments $I_+$ and $I_-$ are identified by an affine map with coefficient $\frac{1}{2i}$. Circular arrows identify adjacent segments.}
\label{fig:dil}
\end{figure}

In Figure 2, the surface is a square with specific identifications of the segments on its horizontal and vertical sides. Circular arrows correspond to conical singularities of angles $\frac{\pi}{2}$ for four of them and of angles $\pi$ for one of them. The conical singularity between segments $I_+$ and $I_-$ has a dilation ratio of $\frac{1}{2}$. The last singularity ($A$ on the figure) has a total angle of $5\pi$ and a dilation ratio of $2$. The sum of the angle defects then implies that the surface has genus zero.\newline 

Geodesic (the dotted line in Figure 2) starting from $M$ with a vertical direction has infinitely many self-intersections and becomes closer and closer to the sides of the square (because of the contraction of the segment $I_-$ to the segment $I_+$). Its $\omega$-limit set is a union of finitely many saddle connections. In particular, this example shows that even if a geodesic has infinitely many self-intersections, its $\omega$-limit set might have empty interior.\newline

The latter example presents a geodesic related to a branched affine structure. However Lemma~\ref{def:geodesic} implies that this curve is also a geodesic for the Fuchsian meromorphic connection induced by the branched affine structure (see Theorem~2).\newline

\section{Flat structures and $k$-differentials}\label{sec:flat}

\subsection{General facts}

For a branched affine surface such that all the linear coefficients of its (affine) transition maps are the $k^{th}$ roots of unity, one gets  a global flat structure. (Observe that in such a case the length of a curve is globally well-defined since the holonomy is unitary). In other words, for any local coordinate $f$, the $k$-differential $(df)^{k}$ is globally well-defined on the whole surface. Therefore, this special kind of a branched affine structure is canonically associated to a $k$-differential.

Locally, the union of all the $k^{th}$ roots of a $k$-differential constitute a multi-valued Abelian differential defining a flat structure. This structure consists of  an atlas of complex charts whose transition maps are given by parallel translations and rotations of order $k$ or its multiple. In particular,  for such a structure, the notions of length, area and direction modulo  rotations of order $k$ are well-defined.
 In the flat structure defined by a $k$-differential, a singularity of order $a > -k$ corresponds to a conical singularity of angle $\frac{(a+k)2\pi}{k}$.

\smallskip 
Flat geometry of $k$-differentials has been studied in e.g. \cite{BCGGM, ShTah, Ta}. In particular, $k$-differentials have  \textit{canonical covers} of order $k$ ramified at their singularities whose orders are not divisible by $k$. The \emph{canonical cover} of a Riemann surface $X$ endowed with a $k$-differential is a Riemann surface $\tilde{X}$ (usually of a higher genus) endowed with a $k$-differential which is a global $k^{th}$ power of a meromorphic $1$-form together with a natural projection $\mu: \tilde X \to X$ sending one $k$-differential to the other. Consequently, up to a rotation of order $k$, there is a well-defined translation structure on $\tilde{X}$. Recall that  a \emph{translation structure} is an atlas of complex charts where all transition maps are translations and all singularities are conical points whose angles are integer multiple of $2\pi$, see \cite{Zo}.\newline

Every geodesic on $X$ lifts to $\tilde{X}$ and thus the dynamics of a geodesic flow in the flat structure of the original $k$-differential splits into the dynamics on the translation surface $\tilde X$ and  a finite monodromy in the fiber. The following theorem classifies the invariant components of translation surfaces for the directional foliation. (This result is a a special case of Proposition 5.5 of \cite{Ta}).

\begin{Theorem}\label{theorem:dichotomy}
Cutting a translation surface without boundary which has a finite total area along its saddle connections in a given direction one decomposes the surface into finitely many connected components of the following two types: 

\noindent
{\rm (i)} \textit{Flat cylinders} whose vertical leaves are periodic geodesics;

\noindent
{\rm (ii)} \textit{Minimal components} for which the noncritical vertical leaves are minimal. The geodesic dynamics on each minimal component is that of an interval exchange map.
\end{Theorem}

The latter dichotomy generalizes the distinction between foliations on flat tori with rational and irrational slopes. As opposed to the hyperbolic cylinders (see Section 3), flat cylinders are obtained by identifying a pair of sides of a parallelogram by a translation.\newline

The only part of Theorem~\ref{theorem:dichotomy} we will need later is the case of a direction in which there is no saddle connection.

\begin{coro+}\label{cor:minimaleaf}
In a translation surface without boundary having a finite area, in a direction without saddle connections, every leaf of the directional foliation is dense in the surface.
\end{coro+}

\begin{proof}
Invariant components in a given direction are bounded by saddle connections. In a direction without any saddle connection, either the whole surface is a minimal component or the whole surface is a flat cylinder. A flat cylinder which is not bounded by any saddle connection is of infinite area.
\end{proof}

Minimality of geodesics in generic directions in the case of translation surfaces extends easily to the case of $k$-differentials.

\begin{coro+}\label{cor:minimalek}
For a Riemann surface $X$ endowed with a $k$-differential inducing a flat structure with finite area, a geodesic in a generic direction  accumulates on the whole surface and its self-intersections are everywhere dense.
\end{coro+}

\begin{proof}
Directions in the affine structure of $X$ are well-defined up to a rotation of order $k$. Every singularity of $X$ is  conical since otherwise the flat structure would have infinite area. Besides, in every homotopy class of topological arcs  between conical singularities (with possible self-intersections), there is at most one saddle connection. Therefore, directions of saddle connections are countable in the circle of directions (modulo rotation of order $k$). These saddle connections lift to saddle connection of the canonical cover $\tilde{X}$ of $X$. Since $\tilde{X}$ is a translation surface, the lift of a geodesic $\gamma$ of $X$ with a generic direction is a geodesic $\tilde{\gamma}$ of $\tilde{X}$ with a direction without any saddle connection. Corollary~\ref{cor:minimaleaf} then implies that $\tilde{\gamma}$ is minimal in $\tilde{X}$. Since $\gamma$ is the projection of $\tilde{\gamma}$, it accumulates everywhere and its self-intersection is also everywhere dense.
\end{proof}

\subsection{Cut-and-paste constructions}

Branched affine structures are well-suited to perform surgeries. Cutting along saddle connections on a surface with a branched affine structure, we obtained a surface with boundary. Moreover, this boundary is geodesic for the underlying affine structure.\newline

For two such affine surfaces with boundary $X$ and $Y$, an identification of corresponding saddle connections by affine maps provides a branched affine structure on $X \cup Y$.\newline

In the case we are interested in, $X$ and $Y$ have affine structures induced by  $k$-differentials (the multiplicators of the monodromy are $k^{th}$ roots of the unity). In this case, if the boundary saddle connections are identified by a composition of translations and rotations of order $k$, then the affine structure on $X \cup Y$ also has a constrained monodromy. %\newline
Since $X \cup Y$ (punctured at the branching points) has a complex affine atlas, it has a natural structure of a Riemann surface. Besides, since the  monodromy is constrained the affine structure  on $X \cup Y$ is induced by some $k$-differential.\newline

Though the surgery is easy to perform and it has a clear geometric meaning, we have no control on the resulting complex structure of $X \cup Y$ (its genus and location of the singularities). In particular, there is no obvious way to give an explicit expression of the $k$-differential.\newline

On a Riemann surface, a branched affine structure (like the one produced by a cut-and-paste surgery) induces a  meromorphic Fuchsian connection (see Theorem~2).

\subsection{First example: quartic differential on $\bC P^1$}

\begin{figure}
\centering
\includegraphics[scale=0.4]{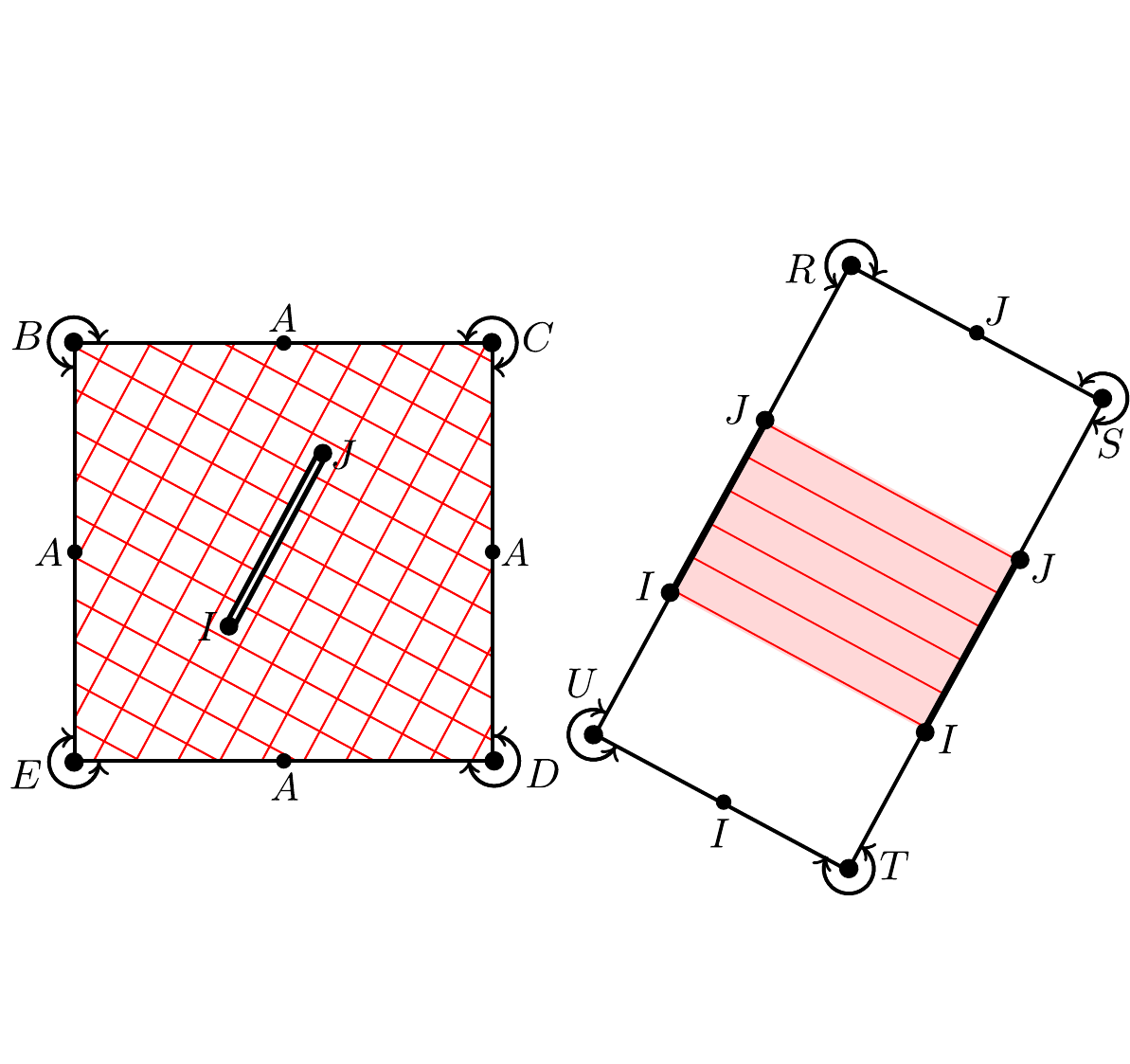} 
%\end{center}

\vskip 0.5cm

\caption {Quartic differential on $\bC P^1$.}
\label{fig:ABquartic}
\end{figure}

The branched affine surface $X$ illustrated in Fig.~\ref{fig:ABquartic} is obtained by a cut-and-paste construction. The left part $X_{l}$ and the right part $X_{r}$ are branched affine surfaces with geodesic boundary.\newline
The left part $X_{l}$ is a square whose horizontal and vertical sides are identified in a specific way (the circular arrows indicate the pairs of identified segments). Corners $B$, $C$, $D$ and $E$ are conical singularities of angle $\frac{\pi}{2}$. The singularity $A$ has a total angle of $4\pi$. The boundary is made of the two sides of a slit $IJ$. Additionally, the slit belongs to a generic direction. Since the identifications are rotations of angle $\frac{\pi}{2}$, the branched affine structure comes from a quartic differential. A computation of the total angle defect shows that the $X_{l}$ jas genus zero and one boundary component.\newline
The right part $X_{r}$ is a rectangle whose sides are identified in a way indicated by the circular arrows. Corners $R$, $S$, $T$, $U$ are conical singularities of angle $\frac{\pi}{2}$. The two lateral segments $IJ$ form the boundary of this branched affine surface. $I$ and $J$ are conical singularities of angle $3\pi$. One can check that $X_{r}$ is also a surface of genus zero with one boundary component. Its branched affine structure is also induced by a quartic differential.\newline
The gluing of the left and right parts along segments $IJ$ preserves lengths and directions. Therefore, the obtained surface $X$ is $\mathbb{CP}^{1}$ with a branched affine structure also induced by a quartic differential.\newline

\begin{prop+}
There is a Fuchsian meromorphic connection on $\mathbb{CP}^{1}$ for which a geodesic $\gamma$ has infinitely many self-intersections and is not everywhere dense. 
\end{prop+}

\begin{proof}
The example of Fig.~\ref{fig:ABquartic} provides a branched affine surface $X$ of genus zero. Its affine structure induces a Fuchsian meromorphic connection (see Theorem~2). We consider a geodesic $\gamma$ starting at the conical singularity $E$ of $X_{l}$ and having the same direction as the slit $IJ$.\newline
The key point is that $\gamma$ can cross the slit $IJ$ and enter $X_{r}$ in only one direction (perpendicular to the slit). Indeed, the flat structure is given by a quartic differential so there are at most two tangent directions at each self-intersection. A branch of $\gamma$ which is parallel to the slit $IJ$ cannot cross it. Then, since $X_{r}$ is a rectangle, if $\gamma$ enters $X_{r}$ through a side of the slit $IJ$, it leaves $X_{r}$ directly  through the other boundary segment (without self-intersections) and crosses the other side of the slit $IJ$. We can glue  portions of $\gamma \cap X_{l}$ to each other and get a geodesic $\gamma'$ of $\bar{X_{l}}$ (a surface without boundary obtained from $X_{r}$ by erasing the slit). By our hypothesis, $\gamma$ never hits a singularity since otherwise it will be a saddle connection which contradicts to the genericity hypothesis on the direction of the slit $IJ$. Following Corollary~\ref{cor:minimalek}, $\gamma'$ is everywhere dense in $\bar{X_{l}}$ and its self-intersections are also everywhere dense. Thus, $\gamma$ accumulates on the union of $X_{l}$ and the rectangle in $X_{r}$. Its self-intersections are dense in $X_{l}$.
\end{proof}

\subsection{Second example: quartic differential on a torus}

\begin{figure}

\begin{center}
\includegraphics[scale=0.55]{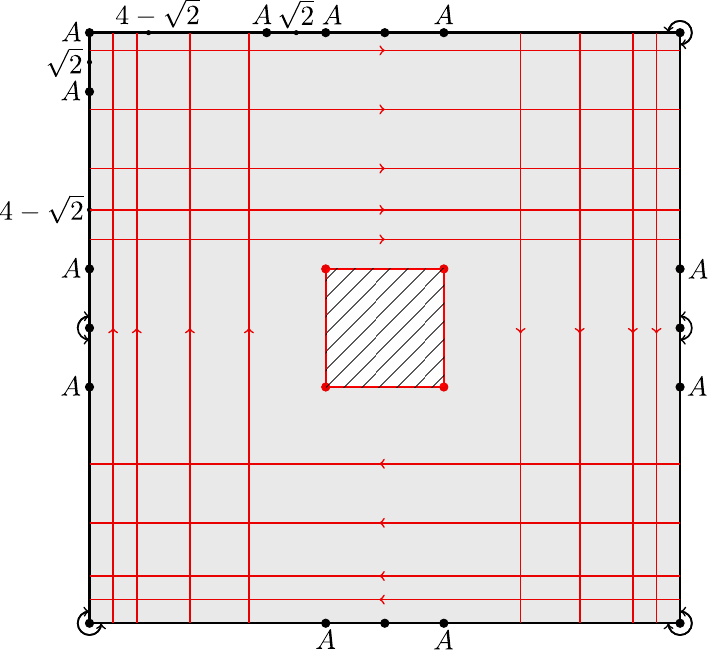} 
\end{center}

\vskip 0.5cm

\caption {Quartic differential on a torus}
\label{fig:AB}
\end{figure}

The branched affine surface $X$ illustrated in Fig~\ref{fig:AB} is a square with boundary segments identified in a special ways. Circular arrows identify the adjacent segments. These identifications create  three conical singularities of angle $\frac{\pi}{2}$ and four conical singularities of angle $\pi$. Besides, two pairs segments of lengths $\sqrt{2}$ and $4-\sqrt{2}$ are identified. The identification maps are compositions of translations and rotations of order four which implies that the branched affine structure is induced by a quartic differential. Singularity $A$ has a total angle of $\frac{21\pi}{2}$. The total angle defect of the singularities on the surface is zero which using the Gauss-Bonnet formula for a flat metric with conical singularities implies that the surface has genus one.\newline
\smallskip 

The canonical cover of this branched affine surface is a translation surface of finite area $\tilde{X}$ (see Subsection 4.1). Horizontal and vertical directions of $X$ are lifted to the same direction on $\tilde{X}$. Theorem~\ref{theorem:dichotomy} then gives its decomposition into invariant components. For this direction, $\tilde{X}$ decomposes into three invariant components two of which being  cylinders. Namely, there is one horizontal and one vertical cylinder in Fig~\ref{fig:AB}. They are formed by periodic geodesics, are bounded by dashed lines and their intersection is the central white square. The third invariant component is a minimal component $\mathcal{M}$ whose dynamics is given by an interval exchange map transposing two segments of lengths $\sqrt{2}$ and $4-\sqrt{2}$.\newline

We consider a geodesic $\gamma$ starting from the point $M$ in the horizontal direction. The lift $\tilde{\gamma}$ of this trajectory belongs to the minimal component $\mathcal{M}$ of $\tilde{X}$. Then, the $\omega$-limit set of $\gamma$ in $X$ is the projection of $\mathcal{M}$ on $X$. This projection is the complement of the small white square at the center of $X$. This central square is the intersection of the projections on $X$ of the two cylinders of $\tilde{X}$. Every horizontal or vertical geodesic of $X$ crossing this square is periodic.\newline

Geodesic $\gamma$ is another example of a geodesic with infinitely many self-intersections whose $\omega$-limit set avoids some part of the surface. In fact, any other geodesic whose lift belongs to the same minimal component could have been chosen.\newline

\section{Final remarks}

All the previous examples can be generalized to surfaces of higher genera by using the following local surgery. Assume that we have a surface with a connection and its geodesic which is not everywhere dense. Take an arbitrarily small disk disjoint from the closure of the latter geodesic and remove a quadrilateral in this disc. Then identify the opposite sides of what was the boundary of the removed quadrilateral. This surgery provides a branched affine surface corresponding to a  meromorphic Fuchsian connection (see Theorem~2 and Subsection 4.2). It increases the genus of the underlying surface by one, but does not impact the given geodesic. Therefore, the behaviour observed in our previous examples exists in all genera.\newline

Besides, as the sides of the quadrilateral can be given arbitrary length and directions, the monodromy of the branched affine structure can be made more generic. Such a surgery proves that examples of Subsections 4.3 and 4.4 are not specific to quartic differentials.\newline

\end{document}